\documentstyle{amsppt}
\magnification=\magstep1
\NoRunningHeads
\NoBlackBoxes


%
%

\pageheight{9truein}
\pagewidth{6.5truein}
\TagsOnRight
\define\lrfloor#1{{\left\lfloor#1\right\rfloor}}
\define\lra#1{{\left\langle#1\right\rangle}}
\define\lrbrace#1{{\left\{#1\right\}}}
\define\xp#1{{\phantom{#1}}}
\define\({\left(}
\define\){\right)}
\define\[{\left[}
\define\]{\right]}
\define\hra{\,{\hookrightarrow}\,}


\let\dsp=\displaystyle
\let\germa=I
\define\vector#1#2{x_{#1},\dots,x_{#2}}
\define\vecy#1#2{y_{#1},\dots,y_{#2}}
%

%
%


%



%

%


\define\ocal{{\Cal O}}

\define\bx{{\bold x}}

\define\ann{\operatorname{ann}}
\define\Ass{\operatorname{Ass}}

\define\pd{\operatorname{pd}}


\define\Coker{\operatorname{Coker}}
\define\depth{\operatorname{depth}}
\define\grade{\operatorname{grade}}
\define\support{\operatorname{support}}
\define\syz{\operatorname{syz}}
\define\Syz{\operatorname{Syz}}
\define\height{\operatorname{height}}

\define\Ext{\operatorname{Ext}}
\define\Hom{\operatorname{Hom}}
\define\im{\operatorname{Im}}
\define\Ker{\operatorname{Ker}}

\define\rank{\operatorname{rank}}
\define\spec{\operatorname{Spec}}

\define\Tor{\operatorname{Tor}}

\topmatter
\title On Modules of Finite Projective Dimension\endtitle
\author S.~P.\ Dutta 
\endauthor
\address\newline
Department of Mathematics\newline
University of Illinois\newline
1409 West Green Street\newline
Urbana, IL 61801\newline
U.S.A.
\endaddress
%
\thanks
\hbox{\hskip-12pt}AMS Subject Classification: Primary 13D02, 13D22;
Secondary 13C15, 13D25, 13H05\hfill\hfill\newline
Key words and phrases: Projective dimension, grade, order ideal,
syzygy, big Cohen-Macaulay module.\hfill\hfill 
\endthanks

\abstract
We address two aspects of finitely generated modules of finite projective dimension
over local rings and their connection in between: embeddability and grade of order ideals
of minimal generators of syzygies. We provide a solution of the embeddability problem
and prove important reductions and special cases of the order ideal conjecture. In particular
we derive that in any local ring $R$ of mixed characteristic $p > 0$, where $p$ is a non-zero-divisor,
if $I$ is an ideal of finite projective dimension over $R$ and $p \ \in \ I$ or $p$ is a
non-zero-divisor on $R/I$, then every minimal generator of $I$ is a non-zero-divisor. Hence if $P$
is a prime ideal of finite projective dimension in a local ring $R$, then every minimal generator
of $P$ is a non-zero-divisor in $R$.
\endabstract

\endtopmatter

\document


\bigskip

\baselineskip17pt


In this note we would like to consider two aspects of finitely
generated modules of finite projective dimension over any local ring:
embeddability and grade of order ideals of minimal generators of its syzygies (minimal). In
regard to embeddability Auslander and Buchweitz ([A-B]) proved that any finitely
generated module on a Gorenstein local ring can be embedded in a module of finite
projective dimension such that the cokernel is Cohen-Macaulay. This result has
several applications in solving homological questions and conjectures in commutative algebra including
Serre's $\chi_i$-conjectures and its generalizations on intersection multiplicity ([D3]). For all
these conjectures one is usually concerned with finitely generated modules $M$ which are of finite projective
dimension over a local ring $R$, but are not so over $R/xR$, $x$ being a non-zero-divisor in
the annihilator of $M$ in $R$. However, a similar result has been absent for non-Gorenstein local rings until now.
In this note we prove the following with respect to embeddability for modules of finite projective dimension:

\proclaim{Theorem\enspace 1.2}Let $(R,m)$ be a local ring and let $M$
be a finitely generated module of finite projective dimension with
positive grade over $R$. Let $\bx=\{x_1,\dots,x_t\}$ be any $R$-sequence
contained in the annihilator of $M$ (henceforth $\ann_RM$) and let
$\overline{R}=R/(x_1,\dots,x_t)$. Then there exists a short exact
sequence of finitely generated $\overline{R}$-modules
$$
0\to M\to Q\to T\to 0
$$
where $\pd_{\overline{R}}Q<\infty$ and $\pd_RT=t$.
\endproclaim

From the construction of $Q$ it would follow that support of $Q$ over
$\overline{R}=\spec(\overline{R})$ and $Q$ satisfies the strong
intersection conjecture due to Peskine and Szpiro (\cite{P-S1}) (Remark 1.2). Since
$T$ is perfect, $T$ also possesses the above property.
The effect of this theorem on generalizations of Serre's conjectures
on intersection multiplicities for arbitrary local rings, in particular
for Cohen-Macaulay rings, would be the subject matter of a future paper.
In this paper we focus on its effect on the order ideal conjecture which
is our next topic.

The order ideal conjecture stems from Evans and Griffith's work on grade of order
ideals of minimal generators of syzygies in equicharacteristic. The statement of the
conjecture is the following.

\proclaim{Order Ideal Conjecture}Let $(R,m)$ be a local ring. Let
$M$ be a finitely generated module of finite projective
dimension over $R$ and let $S_i$ denote its $i^{th}$ syzygy for $i>0$. If $\beta$ is
a minimal generator of $S_i$, then the order ideal $\ocal_{S_i}(\beta)$ has grade at least $i$.
\endproclaim

Let us recall that $\ocal_{S_i}(\beta)=\{f(\beta)|f\in \Hom_R({S_i},R)\}$.

We say that a module $M$ satisfies the order ideal conjecture if order ideals of
minimal generators of all its syzygies satisfy the respective grade inequalities
mentioned above.

Evans and Griffith (\cite{E-G1}, \cite{E-G2}) proved the above
conjecture for equicharacteristic local rings in order to solve the
syzygy problem over the above class of rings. The existence of big
Cohen-Macaulay modules, due to Hochster ([H1]), played an important role in
their proof. Later they proved a graded version of the above
conjecture for a certain class of graded rings in mixed characteristic
(\cite{E-G3}). We would also refer the reader to theorem 9.5.2 in [Br-H]
for a more general version of the order ideal theorem in the equicharacteristic
case. Actually for their proof of the syzygy theorem in equicharacteristic Evans and
Griffith only needed to prove the above conjecture for modules which are locally
free on the punctured spectrum of spec(R), R being regular local. And they reduced the proof
of this case to what is now known as improved new intersection conjecture. Later
Hochster showed that the canonical element conjecture implies the improved new intersection
conjecture [H2]. The equivalence of these two conjectures was established in [D1].
In our most recent work [D2] we have shown that \it a particular case of the order ideal conjecture
implies the monomial conjecture \rm and hence all its equivalent forms, e.g., the direct
summand conjecture, the canonical element conjecture, the improved intersection conjecture etc.
Thus the order ideal conjecture now occupies a central
position among several homological conjectures in commutative algebra.

First let us mention that in order to prove the order ideal conjecture
on arbitrary local rings $R$, \it it is enough to concentrate on first syzygies \rm of
modules of finite projective dimension (Lemma 2.1).  Theorem 2.3 shows that \it for the
validity of the order ideal conjecture it is enough to prove that every minimal generator
of ideals of height 2, grade 2 and of finite projective dimension
over $R$ is a non-zero-divisor. \rm In theorem 2.5 we prove the following:

\it
Given a module of finite projective dimension $M$ on
a local ring $R$, there exists an $R$-sequence
$x_1,\dots,x_h$, $h=\pd_RM$, such that for any $R$-module $N$ for
which $x_1,\dots,x_h$ form an $N$-sequence, $\Tor_j^R(M,N)=0$ for
$j>0$.

\rm
This theorem leads us to the following (Corollary 2, 2.7):

\it
Let $(R,m)$ be a local ring and $M$ be a
finitely generated module of finite projective dimension over $R$. Let
$x_1,\dots,x_h$ be an $R$-sequence as mentioned in the above
proposition. If $R$ has mixed characteristic $p>0$, we assume that
$p,\,x_1,\dots,x_h$ form a part of a system of parameters of $R$. Then
for every minimal generator $\beta$ of $S_1=\Syz_R^1(M)$, $\grade
\ocal_{S_1}(\beta)\ge 1$.
\medskip
\rm
The statement of our main theorem in section 2 is the following.

\proclaim{Theorem\enspace 2.11}Let $(R,m)$ be a local ring of mixed
characteristic $p>0$. Let $M$ be a finitely generated module of finite
projective dimension over $R$ and let $\beta$ be a minimal generator
of $S_i$, the $i^{\text{th}}$ syzygy (minimal) of $M$, for $i>0$. We
assume that either $p$ is a non-zero-divisor in $R$ or $p$ is
nilpotent. We have the following:

\roster
\item"a)" if $p$ is nilpotent, then the order ideal conjecture is
valid on~$R$.

\item"b)" if $pM=0$, then grade of $\ocal_{S_1}(\beta)\ge 1$.

\item"c)" Suppose that balanced or complete almost Cohen-Macaulay algebras
exist over complete local domains. If $pM=0$ and $p\in m-m^2$, then grade of
$\ocal_{S_i}(\beta)\ge i$, $\forall i \ge 1$.

\item"d)" Assume that every element in \ $m-m^2$ is a non-zero-divisor and the order
ideal conjecture is valid over $R/xR$ for $x  \in  m-m^2$. If $\ann_RM \cap (m-m^2) \ne \phi$
or $depth_R M > 0$, then $grade\ocal_{S_i}(\beta)\ge i$, $\forall i \ge 1$.
\endroster
\endproclaim

As a consequence of the above theorem we have the following corollaries:

\proclaim{Corollary\enspace 2.12} For any ideal $I$ of finite projective dimension over $R$ of mixed characteristic
$p > 0$, where $p$ is a non-zero-divisor in $R$, if $p  \in  I$ or $p$ is a
non-zero-divisor on $R/I$ then every minimal generator of $I$ is a non-zero-divisor.
In particular if $P$ is a prime ideal of finite projective dimension
over $R$, then every minimal generator of $P$ is a non-zero-divisor in $R$.
\endproclaim

\proclaim {Corollary\enspace 2.13}\enspace Let (R, m) be a regular local ring of dimension $n$ and assume that
the order ideal conjecture is valid for regular local rings of dimension $(n-1)$. If $M$ is a finitely
generated $R$-module such that either $M$ is annihilated by a regular parameter or $\depth_R M > 0$,
then $M$ satisfies the order ideal conjecture.
\endproclaim

The two main ingredients of our proof of theorem (2.11) are theorem 1.2 and
Shimomoto's theorem (2.10) on existence of almost Cohen-Macaulay algebras (\cite{Shi}).



Throughout this note "local" means noetherian local. For definitions of standard
notions like projective dimension, grade etc. and their basic properties we refer the reader to [Br-H].

\head{Section 1}\endhead

First we would like to mention the following proposition.

\proclaim{1.1\enspace Proposition\enspace {\rm (Prop.\ 1.1,
\cite{D1})}} Let $A$ be a noetherian local ring
and let $F_{\bullet}:\to A^{s_i}\to A^{s_{i-1}}\to \cdots \to
A^{s_0}\to 0$ be a free complex with $H_0(F_{\bullet}) = M$. Let
$N\subset M$ be a submodule of $M$. Then we can construct a free
complex $G_{\bullet}:\to A^{d_i}\to A^{d_{i-1}}\to \cdots \to
A^{d_0}\to 0$ with $H_0(G) = N$ and a map
$\phi_{\bullet}:G_{\bullet}\to F_{\bullet}$ such that $\phi_0$ induces
the inclusion $N\hra M$ and the mapping cone of $\phi_{\bullet}$ is a
free resolution of $M/N$ i.e., $\phi_{\bullet}$ induces an isomorphism
$H_i(G_{\bullet})\widetilde{\to} H_i(F_{\bullet})$ for
$i>0$. Moreover, by our construction, $G_{\bullet}$ is minimal.
\endproclaim

For a proof we refer the reader to Prop.\ 1.1 in \cite{D1}.

Next we prove the main theorem of this section.

\proclaim{1.2.\enspace Theorem}Let $\(R, m\)$ be a local ring and let $M$ be a
finitely generated module of finite projective dimension over $R$. Assume that
$\grade_RM>0$. Let ${\bold{x}}$ denote the ideal generated by
an $R$-sequence of length $t$ contained in $\ann_R M$ and
let $\overline{R} = R/{\bold{x}}R$. Then there exists a short exact
sequence of finitely generated $\overline{R}$-modules
$$
0\to M\to Q\to T\to 0
$$
such that $\pd_{\overline{R}}Q<\infty$ and $\pd_RT=t$.
\endproclaim

\demo{Proof}First suppose that $\pd_R M = 1$. Let us recall that $\grade_RM \leq \pd_RM$.
Since $\grade_R M > 0$, it follows that $\grade_RM = 1$. Let $x \in \ann_R M$ be a
non-zero-divisor. Consider a minimal free resolution: $0 \to R^{t_1} @>{f}>>\,R^{t_0}
\to M \to 0$ \  of \ M over R. Tensoring this resolution with $R/xR$ we obtain an exact
sequence: $0 \to Tor_1^R (M, R/xR) \to (R/xR)^{t_1} @>{\overline{f}}>>\, (R/xR)^{t_0}
\to M \to 0$. Since $x \in \ann_R M, Tor_1^R (M, R/xR) \simeq M$. Let $T = im \overline{f}$.
Hence we obtain the following short exact sequence: $0 \to M \to (R/xR)^{t_1} \to T \to 0$
and our assertion follows. So we can assume that $\pd_R M \ge 2$.

Let $(F_{\bullet}, d_{\bullet}):0\to R^{r_n}
@>{d_n}>>\,R^{r_{n-1}}\to\cdots\to R^{r_1} @>{d_1}>> R^{r_0}\to 0$ be a
minimal projective resolution of $M$ over $R$ and let
$(L_{\bullet}, \phi_{\bullet}):\to \overline{R}^{s_n} @>\phi_n>>
\overline{R}^{s_{n-1}}\to \cdots \overline{R}^{s_1} @>\phi_0>>
\overline{R}^{s_0}\to 0$ be a minimal projective resolution of $M$
over $\overline{R}$. Since $\pd_RM=n$, $\grade \Ext^i_R(M,R)\ge i$
for $1\le i\le n$, $\Ext^n_R(M,R)\ne 0$ and $\Ext^i_R(M,R)=0$ for
$i>n$. Hence $\grade \Ext^i_{\overline{R}}(M,\overline{R})\ge i$ for
$1\le i\le n-t$, $\Ext^{n-t}_{\overline{R}}(M,\overline{R})\ne 0$
$(\simeq \Ext^n_R(M,R))$ and $\Ext^i_{\overline{R}}(M,\overline{R})=0$
for $i>(n-t)$.

Applying $\Hom_{\overline{R}}(-,\overline{R})$ to $L_{\bullet}$ we
obtain the following free complex $L^\ast_{\bullet}$:
$$
L^\ast_{\bullet}:0\to \overline{R}^{s_0^\ast}\to
\overline{R}^{s_1^\ast}\to\cdots\to \overline{R}^{s^\ast_{n-t}}\to
0.
$$
(For any $R(\overline{R})$ module $V$,
$V^\ast=\Hom_R(M,R)(\Hom_{\overline{R}}(V,\overline{R}))$)

Let $G=\Coker \phi^\ast_{n-t}$. We have a short exact sequence
$$
0\to \Ext^{n-t}_{\overline{R}}(M,\overline{R}) @>c>> G @>\eta>> \im
\phi^\ast_{n-t+1}\to 0.
\tag{1}
$$
By the above proposition, there exists a minimal free complex
$(P_{\bullet},\alpha_{\bullet}):\to \overline{R}^{b_{n-t}}
@>\alpha_{n-t}>> \overline{R}^{b_{n-t-1}}\to\cdots\cdots @>\alpha_1>>
\overline{R}^{b_0}\to 0$ with
$H_0(P_{\bullet})=\Ext^{n-t}_{\overline{R}}(M,\overline{R})$ and a map
$\Psi_{\bullet}:P_{\bullet}\to L_{\bullet}^\ast$ such that
$\Psi_{\bullet}$ induces the injection $c$ in (1) and the mapping cone
of $\Psi_{\bullet}$ is a free resolution of $\im \phi_{n-t+1}^\ast$. We have
the following commutative diagram


$$
\vbox{\hbox{\vbox{\halign{\hbox to 28pt{}\eightpoint\tabskip=3pt\!
$#$\hfill\hfill & $#$\hfill\hfill & $#$\hfill\hfill & $#$\hfill\hfill &
$#$\hfill\hfill & $#$\hfill\hfill & $#$\hfill\hfill & $#$\hfill\hfill &
$#$\hfill\hfill & $#$\hfill\hfill & $#$\hfill\hfill & $#$\hfill\hfill &
$#$\hfill\hfill & $#$\hfill\hfill & $#$\hfill\hfill & $#$\hfill\hfill &
$#$\hfill\hfill & $#$\hfill\hfill\tabskip=0pt\cr
P_\bullet: & \longrightarrow & \overline{R}^{b_{n-t+1}} & \hbox{$\dsp
@>\,\alpha_{n-t+1}\,>>$} & \overline{R}^{b_{n-t}} & \longrightarrow &
\overline{R}^{b_{n-t-1}} & \cdots & \longrightarrow & \overline{R}^{b_1} &
\longrightarrow & \overline{R}^{b_0} & \longrightarrow & 0\cr
&&&& \Big\downarrow{\Psi_{n-t}} && \Big\downarrow{\Psi_{n-t-1}} &&&
\Big\downarrow{\Psi_1} && \Big\downarrow{\Psi_0}\cr
L_\bullet^\ast:&& 0 & \longrightarrow & \overline{R}^{s^\ast_0} &
\longrightarrow & \overline{R}^{s_1^\ast} & \cdots & \longrightarrow &
\overline{R}^{s^\ast_{n-t-1}} & @>{\phi^\ast_{n-t}}>> &
\overline{R}^{s^\ast_{n-t}} & \longrightarrow & 0\cr
\noalign{\vskip-34pt}
\noalign{\hfill\hfill\text{(2)}}\cr
\noalign{\vskip10pt}
}}}}
$$
\tenpoint\baselineskip17pt

Applying $\Hom(-,\overline{R})$ to (2) we obtain the following
commutative diagram of exact complexes
$$
\vbox{\hbox{\vbox{\halign{\hbox to 28pt{}\eightpoint\tabskip=3pt\!
$#$\hfill\hfill & $#$\hfill\hfill & $#$\hfill\hfill & $#$\hfill\hfill &
$#$\hfill\hfill & $#$\hfill\hfill & $#$\hfill\hfill & $#$\hfill\hfill &
$#$\hfill\hfill & $#$\hfill\hfill & $#$\hfill\hfill & $#$\hfill\hfill &
$#$\hfill\hfill & $#$\hfill\hfill & $#$\hfill\hfill & $#$\hfill\hfill &
$#$\hfill\hfill & $#$\hfill\hfill\tabskip=0pt\cr
& \longrightarrow & \overline{R}^{s_{n-t}} & \longrightarrow &
\overline{R}^{s_{n-t-1}} & \longrightarrow & \cdots & \longrightarrow &
\overline{R}^{s_1} & \longrightarrow & \overline{R}^{s_0} &
\longrightarrow & M &\longrightarrow & 0\cr
&& \Big\downarrow{\Psi_0^\ast} && \Big\downarrow{\Psi_1^\ast} &&&&
\Big\downarrow{\Psi^\ast_{n-t-1}} && \Big\downarrow{\Psi^\ast_{n-t}} &&
\Big\downarrow{\Psi}\cr
0 & \longrightarrow & \overline{R}^{b_0^\ast} & \longrightarrow &
\overline{R}^{b_1^\ast} & \longrightarrow & \cdots &
\longrightarrow & \overline{R}^{b^\ast_{n-t-1}} & \hbox{$\dsp
@>\,\alpha^\ast_{n-t}\,>>$} & \overline{R}^{b^\ast_{n-t}} &
\longrightarrow & Q & \longrightarrow & 0\cr
\noalign{\vskip-34pt}
\noalign{\hfill\hfill\text{(3)}}\cr
\noalign{\vskip10pt}
}}}}
$$
\tenpoint\baselineskip17pt
where $Q=\Coker \alpha^\ast_{n-t}$ and $\Psi:M\to Q$ is induced by
$\Psi^\ast_{n-t}$. Since $\grade
\Ext^i_{\overline{R}}(M,\overline{R})\ge i$, $1\le i\le n-t$,
the bottom row of (3) provides a minimal free resolution of
$Q$ and hence $\pd_{\overline{R}}Q=n-t$. Moreover, by
construction, $\Ext^i_{\overline{R}}(Q,\overline{R}) @>\sim>>
\Ext^i_{\overline{R}}(M,\overline{R})$ for $i>0$.

Next we want to prove that $\Psi:M\to Q$ in (3) is injective.

Let $\overline{F}_{\bullet} = F_{\bullet}\otimes \overline{R}$,
$S_i=\syz^i_R(M)$ and $\overline{S}_i=S_i/{\bold{x}}S_i$. Since
${\bold{x}}$ is generated by an $R$-sequence of length $t$ in $\ann_RM$, we have
$\Tor^R_t(M,\overline{R})\,{\overset{\sim}\to{\longleftarrow}}\,M$ and
$\Tor^R_i(M,\overline{R})=0$ for $i>t$. Tensoring the exact sequence
 $0\to S_t\to R^{r_{t-1}}\to
S_{t-1}\to 0$ by $\overline{R}$ we obtain an exact sequence
$$
0\to M@>j>> \overline{S}_t\to \overline{R}^{r_{t-1}}\to \overline{S}_{t-1}\to 0.
$$
Let $\theta_{\bullet}:L_{\bullet}\to \overline{F}_{\bullet}$ be a lift
of $j:M\hra \overline{S}_t$. We have the following commutative diagram:
$$
\vbox{\hbox{\vbox{\halign{\hbox to 48pt{}\eightpoint\tabskip=3pt\!
$#$\hfill\hfill & $#$\hfill\hfill & $#$\hfill\hfill & $#$\hfill\hfill &
$#$\hfill\hfill & $#$\hfill\hfill & $#$\hfill\hfill & $#$\hfill\hfill &
$#$\hfill\hfill & $#$\hfill\hfill & $#$\hfill\hfill & $#$\hfill\hfill &
$#$\hfill\hfill & $#$\hfill\hfill & $#$\hfill\hfill & $#$\hfill\hfill &
$#$\hfill\hfill & $#$\hfill\hfill\tabskip=0pt\cr
L_\bullet: &&& \longrightarrow & \overline{R}^{s_{n-t}} & \longrightarrow &
\overline{R}^{s_{n-t-1}} & \longrightarrow & \cdots & \longrightarrow &
\overline{R}^{s_0} & \longrightarrow & 0\cr
&&&& \Big\downarrow{\theta_{n-t}} && \Big\downarrow{\theta_{n-t-1}} &&&&
\Big\downarrow{\theta_0}\cr
\overline{F}_\bullet: && 0 & \longrightarrow & \overline{R}^{r_n} &
\longrightarrow & \overline{R}^{r_{n-1}} & \longrightarrow & \cdots &
\longrightarrow & \overline{R}^{r_t} & \longrightarrow &
\overline{R}^{r_{t-1}} & \longrightarrow & \dots\cr
\noalign{\vskip-34pt}
\noalign{\hfill\hfill\text{(4)}}\cr
\noalign{\vskip10pt}
}}}}
$$
\tenpoint\baselineskip17pt

The mapping cone of $\{\theta_\bullet\}$ in (4) is a free resolution
of $\im\,\overline{d}_t$ over $\overline{R}$ and $\theta_0$ induces
the isomorphism $\tilde{\theta}_0:M@>\sim>> \Tor^R_t(M,\overline{R})$
via $j:M\hra \overline{S}_t$.

Applying $\Hom(-,\overline{R})$ to (4) the following commutative
diagram is obtained:
$$
\vbox{\hbox{\vbox{\halign{\hbox to 18pt{}\eightpoint\tabskip=3pt\!
$#$\hfill\hfill & $#$\hfill\hfill & $#$\hfill\hfill & $#$\hfill\hfill &
$#$\hfill\hfill & $#$\hfill\hfill & $#$\hfill\hfill & $#$\hfill\hfill &
$#$\hfill\hfill & $#$\hfill\hfill & $#$\hfill\hfill & $#$\hfill\hfill &
$#$\hfill\hfill & $#$\hfill\hfill & $#$\hfill\hfill & $#$\hfill\hfill &
$#$\hfill\hfill & $#$\hfill\hfill\tabskip=0pt\cr
(\overline{F}_\bullet^\ast) & : & & \longrightarrow &
\overline{R}^{r^\ast_{t-1}} & \hbox{$\dsp @>\,\overline{d}_t^\ast\,>>$} &
\overline{R}^{r^\ast_t} & \hbox{$\dsp
@>\,\overline{d}^\ast_{t+1}\,>>$} & \cdots &
\longrightarrow & \overline{R}^{r^\ast_{n-1}} & \hbox{$\dsp
@>\,\overline{d}^\ast_n\,>>$} & \overline{R}^{r^\ast_n} & \longrightarrow &
\Ext^n_R(M,\overline{R}) & \longrightarrow & 0\cr
&&&&&& \Big\downarrow{\theta_0^\ast} &&&& \Big\downarrow{\theta^\ast_{n-t-1}} &&
\Big\downarrow{\theta^\ast_{n-t}} &&
\Big\downarrow^{\hbox{\hskip-4.25pt}^{_\cap}}\cr
(L_\bullet^\ast) & : & & & 0 & \longrightarrow & \overline{R}^{s^\ast_0} &
\longrightarrow & \cdots & \longrightarrow & \overline{R}^{s^\ast_{n-t-1}} &
\longrightarrow & \overline{R}^{s^\ast_{n-t}} & \longrightarrow & G
\longrightarrow 0\cr
\noalign{\vskip-34pt}
\noalign{\hfill\hfill\text{(5)}}\cr
\noalign{\vskip10pt}
}}}}
$$
\tenpoint\baselineskip17pt

Since the mapping cone of (4) is a free resolution of
$\im\,\overline{d}_t$ over $\overline{R}$ and
$\Ext^i_{\overline{R}}(M,\,\overline{R})=0$ for $i>n-t$,
$\theta^\ast_{n-t}$ induces an isomorphism
$\Ext^n_R(M,\overline{R})@>\sim>>
\Ext_{\overline{R}}^{n-t}(M,\overline{R})$. Hence the inclusion
$\Ext^n_R(M,\overline{R})\hra G$ in (5) can be identified with
$c:\Ext^{n-t}_{\overline{R}}(M,\overline{R})\hra G$ in~(1).

Let $(H_{\bullet},\delta_{\bullet})$ denote the mapping cone of
$\Psi_{\bullet}$ in (2). By construction this is a free resolution of
$\im\phi^\ast_{n-t+1}$. Let $\eta_{\bullet}:L^\ast_{\bullet}\to
H_{\bullet} (\eta_i:\overline{R}^{s^\ast_i}\to
\overline{R}^{b_{n-i-t-1}}\oplus \overline{R}^{s^\ast_i})$,
$\gamma_{\bullet}:H_{\bullet}\to P_{\bullet}(-1)
(\gamma_{n-j}:\overline{R}^{b_{n-j}}\oplus
\overline{R}^{s^\ast_{j-t-1}}\to \overline{R}^{b_{n-j}})$ denote the
corresponding inclusion and projection maps respectively. Then
$\eta_{\bullet}\cdot \theta_{\bullet}^\ast: \overline{F}_{\bullet}^\ast\to
H_{\bullet}$ lifts the composite
$\eta\cdot c:\Ext^{n-t}_{\overline{R}}(M,\overline{R})@>c>> G @>\eta>>
\im\phi^\ast_{n-t+1}$. Since $\eta\cdot c=0$ and
$(H_{\bullet},\delta_{\bullet})$ is a free resolution of $\im
\phi_{n-t+1}^\ast$, $\eta_{\bullet}\cdot \theta_{\bullet}^\ast$ is homotopic
to 0. Hence there exists homotopy maps
$h_{\bullet}:\overline{F}_{\bullet}^\ast \to H_{\bullet}$,
$h_j:\overline{R}^{r_j^\ast}\to \overline{R}^{b_{n-j}}\oplus
\overline{R}^{s^\ast_{j-t-1}}$, for $t\le j\le n$, such that
$$
\delta_{n-j}\cdot h_j+h_{j+1}\cdot \overline{d}^\ast_{j+1} =
\eta_{j-1}\cdot \theta^\ast_{j-1}.
\tag{6}
$$
Let $\beta_{\bullet} = \gamma_{\bullet}\cdot h_{\bullet}$,
$\left\{\beta_j=\gamma_{n-j}\cdot h_j\right\}$. Then
$\beta^\prime_\bullet=\{(-1)^{j+1}\beta_j\}:
\overline{F}_{\bullet}^\ast \to P_\bullet$ is a map of complexes. Consider
$\Psi_{\bullet}.\beta^\prime_{\bullet}:F^\ast_{\bullet}\to L^\ast_{\bullet}$.

For $1\le j\le n$, let $k_j=\pi_2\cdot h_j:\overline{R}^{r_j^\ast}\to
\overline{R}^{s^\ast_{j-t-1}}$, where
$\pi_2:\overline{R}^{b_{n-j}}\oplus \overline{R}^{s^\ast_{j-t-1}}\to
\overline{R}^{s^\ast_{j-t-1}}$ is the projection on the second
component. It can be checked that, for $t\le j\le n$,
$$
\theta^\ast_{j-1}-(-1)^{n-j-1}\Psi_{n-j}\cdot\beta_j=\phi_{j-1}^\ast
\cdot k_j + k_{j+1}\cdot \overline{d}^\ast_{j+1};\tag{7}
$$
i.e.\ $\theta^\ast_\bullet$ and $\Psi_\bullet \beta_\bullet^\prime$
are homotopic.

\baselineskip18pt

The commutative diagrams below are provided to clarify (6) and (7) for
$t=1$.
$$
\vbox{\hbox{\vbox{\halign{\eightpoint\tabskip=2pt\!
$#$\hfill\hfill & \hfill\hfill $#$ & \hfill$#$\hfill &
$#$\hfill\hfill\tabskip=0pt\cr
\overline{F}_\bullet^\ast \hfill\hfill : & \cdots \longrightarrow
\;\;\overline{R}^{r_j^\ast}\, \hbox{$\dsp
 @>\overline{d}^\ast_{j+1}>>$}\;\; \overline{R}^{r^\ast_{j+1}}
 \longrightarrow\quad\;\; & \cdots & \longrightarrow
 \overline{R}^{r^\ast_{n-1}} \longrightarrow \overline{R}^{r^\ast_n}
 \longrightarrow 0\cr
%
%
\noalign{$\dsp \hbox to
116pt{}{}^{hj}\hbox{\hskip-5pt}\Bigg/\hbox to
5pt{}\Bigg\downarrow\hbox{\hskip-2pt}\theta^\ast_{j-1}
\hbox{\hskip-2pt}{}^{h_{j+1}}\hbox{\hskip-5pt}\Bigg/\Bigg\downarrow
\hbox to 91pt{}\Bigg\downarrow \hbox to
18pt{}{}^{hn}\hbox{\hskip-5pt}\Bigg/\hbox to 6pt{}\Bigg\downarrow$}\cr
\noalign{\vskip-12pt}
\noalign{$\dsp \hbox to 110.75pt{}\Bigg/\hbox to
39.5pt{}\Bigg/\hbox to 127.5pt{}\Bigg/$}\cr
\noalign{\vskip-54pt}\cr
%
%
L_\bullet^\ast \hfill\hfill : & \cdots \longrightarrow
\overline{R}^{s^\ast_{j-2}} \longrightarrow
\overline{R}^{s^\ast_{j-1}}\, \longrightarrow\; \overline{R}^{s^\ast_j}
\;\;\;\longrightarrow\quad\; & \cdots &
\longrightarrow \overline{R}^{s^\ast_{n-t-1}} \longrightarrow
\overline{R}^{s^\ast_{n-1}} \longrightarrow 0\cr
%
\noalign{$\dsp \hbox to 70pt{}\eta_{j-2}\Bigg\downarrow \hbox to
40pt{}\Bigg\downarrow\eta_{j-1} \hbox to 119pt{} \Bigg\downarrow \hbox to
43pt{}\Bigg\downarrow$}\cr
\noalign{\vskip-9pt}%
%
H_\bullet \hfill\hfill : & \cdots \longrightarrow
\overline{R}^{b_{n-j}}\oplus \overline{R}^{s^\ast_{j-2}}
{\underset{\delta_{n-j}}\to\longrightarrow}\;
\overline{R}^{b_{n-j-1}}\oplus\overline{R}^{s^\ast_{j-1}}
\longrightarrow \cdots & \longrightarrow &
\overline{R}^{b_0}\hbox{\hskip-2pt}\oplus
\hbox{\hskip-2pt}\overline{R}^{s^\ast_{n-t-1}}
\hbox{\hskip-2pt}\longrightarrow\hbox{\hskip-2pt}
\overline{R}^{s^\ast_{n-1}}\cr
\noalign{\vskip-6pt}%
\noalign{$\dsp \hbox to 62pt{}\gamma_{n-j}\Bigg\downarrow \hbox to
84pt{}\Bigg\downarrow\gamma_{n-j-1} \hbox to 71pt{}\Bigg\downarrow$}\cr
\noalign{\vskip-9pt}%
P_\bullet(-1) \hfill\hfill : & \cdots \longrightarrow
\overline{R}^{b_{n-j}}\qquad\hbox{$\dsp
@>\,\alpha_{n-j}\,>>$}
\qquad \overline{R}^{b_{n-j-1}}\;\; \longrightarrow\quad & \cdots
&\longrightarrow \overline{R}^{b_0} \longrightarrow 0\cr
\noalign{\vskip4pt}%
& \delta_{n-j} = \(\alpha_{n-j},\,(-1)^{n-j-1} \Psi_{n-j} +
\phi^\ast_{j-1}\)\qquad\hfill\hfill\cr
}}}}
$$
\tenpoint\baselineskip18pt
\vglue-0.15truein
$$
\vbox{\hbox{\hbox to 110pt{}\vbox{\halign{\eightpoint\tabskip=2pt\!
$#$\hfill\hfill & \hfill\hfill $#$ & \hfill$#$\hfill &
$#$\hfill\hfill\tabskip=0pt\cr
\overline{F}_\bullet^\ast \hfill\hfill : &
\hbox{$\dsp\cdots\;\;@>{\;\quad\;}>>\;\;\overline{R}^{r_j^\ast}\;\;
@>{\;\quad\;}>>\;\;\overline{R}^{r^\ast_{j+1}}\qquad\;\quad$}\cr
\noalign{$\dsp\hbox to 69.75pt{}{}^{k_j}\hbox{\hskip-5pt}\Bigg/\hbox to
0pt{}\Bigg\downarrow\hbox{\hskip-2pt}\beta_j
\hbox{\hskip14.0pt}{}^{k_{j+1}}\hbox{\hskip-5pt}\Bigg/
\hbox to 0pt{}\Bigg\downarrow\hbox{\hskip-2pt}\beta_{j+1}$}\cr
\noalign{\vskip-8pt}
\noalign{$\dsp \hbox to 62.7pt{}\Bigg/\hbox to 40.7pt{}\Bigg/$}\cr
\noalign{\vskip-56pt}\cr
\noalign{\hbox{\hskip89.5pt$\bigg|$\hskip52.0pt$\bigg|$}}\cr
\noalign{\vskip-10pt}
\noalign{\hbox{\hskip71pt${}_{\theta^\ast_{j-1}}
\Bigg\downarrow$\hskip40pt${}_{\theta_j^\ast}\Bigg\downarrow\hbox{\hskip-2pt}$}}\cr
\noalign{\vskip-58pt}
P_\bullet\hfill\hfill : & \dsp\cdots\hbox{\hskip-2pt}
\longrightarrow\hbox{\hskip-2pt}
\overline{R}^{b_{n-j}}\hbox{\hskip8pt}
\hbox{\hskip-3pt$\longrightarrow$\hskip6pt}
\overline{R}^{b_{n-j-1}}\qquad\;\cr
\noalign{\vskip5pt}
\noalign{$\dsp\hbox to
85pt{}\Bigg\downarrow\hbox{\hskip-2pt}\Psi_{n-j} \hbox to
26pt{}\Bigg\downarrow\Psi_{n-j-1}$}\cr
\noalign{\vskip-10pt}
L^\ast_\bullet \hfill\hfill : & @>{\;\quad\;}>>
\overline{R}^{s^\ast_{j-2}}\;@>{\;\quad\;}>>
\hbox{\hskip2pt$\overline{R}^{s^\ast_{j-1}}$} @>{\;\quad\;}>>
\overline{R}^{s^\ast_j}\hbox to 42pt{}\cr
}}}}
$$
\tenpoint\baselineskip17pt


\noindent
Hence $\{\theta_\bullet^\ast\}$ and
$\{\Psi_\bullet\beta^\prime_\bullet\}$ define identical maps (modulo
$\pm$ sign) on homologies of $F_\bullet^\ast$ and $L_\bullet$ and
consequently for the homologies of their respective duals. If
$\tilde\Psi:M\to H^{n-t}(P^\ast_\bullet)$,
$\tilde\beta:H^{n-t}(P_\bullet^\ast)\to \Tor_t^R(M,\overline{R})$ and
$\tilde\theta_0:M@>\sim>> \Tor_t^R(M,\overline{R})$ denote the maps
induced by $\Psi^\ast_{n-t}$, $\beta^\ast_t$ and
$(\theta^\ast_0)^\ast=\theta_0$ respectively, then $\tilde\theta_0 =
\tilde\beta\cdot\tilde\Psi$ (modulo a sign). Since $\tilde\theta_0$ is an isomorphism,
$\Psi =$ (inclusion of $H^{n-t}(P_\bullet^\ast)\hra Q)\cdot \tilde\Psi:M\hra Q$ is injective.

Let $\Coker\psi = T$; then $0\to M @>\Psi>> Q\to T\to$ is exact. Since
the mapping cones of $\Psi_{\bullet}^\ast,\Psi_{\bullet}$ are free
resolutions of $T$, $\im \phi_{n-t+1}^\ast$ respectively and
$\Ext^i_{\overline{R}}(M,\overline{R})@>\sim>>
\Ext^i_{\overline{R}}(Q,\overline{R})$ for $i>0$, it follows that
$\Ext^i_{\overline{R}}(T,\overline{R})=0$ for $i>0$ i.e.,
$\Ext^i_R(T,R)=0$ for $i>t$. Since $\pd_RM<\infty$, we have
$\pd_RT<\infty$ and hence $\pd_RT=t$. This completes our proof.
\enddemo

\proclaim{1.3 \enspace Corollary}Let $R$ be a local ring and let $M$ be a finitely
generated module of finite projective dimension with $\pd_RM>1$ and
$\grade R>0$. Given any non-zero-divisor $x\in \ann_RM$, $M$ can be
imbedded in a finitely generated module $Q$ of finite projective
dimension over $R/xR$ in such a way that if $(G_{\bullet},\gamma_{\bullet})$,
$(F_{\bullet},d_{\bullet})$ are minimal free resolutions of $M$ and
$Q$ over $R$ respectively and $\phi_{\bullet}:G_{\bullet}\to
F_{\bullet}$ is a lift of $i:M\hra Q$, then $\phi_{\bullet}$ induces
our isomorphism between $(G_{\bullet},\gamma_{\bullet})_{i\ge 2}$ and
$(F_{\bullet},d_{\bullet})_{i\ge 2}$, $\phi_1(G_1)$ is either a
summand of or isomorphic to $F_1$ and $\phi_0$ is an
injection. Moreover, $\syz^1_R(M)\oplus R^t\simeq \syz^1_R(Q)$ for
some $t\ge 0$.
\endproclaim

\demo{Proof}By the above theorem we have an exact sequence of $R/xR$ modules
$$
0\to M @>i>> Q @>\eta>> T\to 0
$$
where $\pd_{R/xR}Q<\infty$, $\pd_RT=1$.

The proof of the corollary now follows directly either by constructing
a minimal free resolution of $Q$ from minimal free resolutions of $M$
and $T$ or by extracting a minimal free resolution of $M$ from the
mapping cone of $\eta_{\bullet}:F_{\bullet}\to L_{\bullet}$, where
$F_{\bullet}, L_{\bullet}$ are minimal free resolutions of $Q$ and $T$
respectively. For details of such basic constructions the reader is referred
to [Br-H].
\enddemo

\medskip

\noindent {\bf Remark.}\enspace Since $\grade_R T =  \pd_R T =  r$, $T$ is perfect.
It can be easily checked from the above construction that
$\grade_{\overline R} Q =  0$ and since $\pd_{\overline R} Q < \infty$,
 $\text{support\,}_{\overline{R}}Q = \spec(\overline{R})$. Moreover, the strong
intersection conjecture ([H1], [P-S]) is valid for both $Q$ and~$T$. For details on
this observations we refer the the reader to section 4, chapter II in [P-S].

\head{Section 2}\endhead

\proclaim{2.1\enspace Lemma } Let (R, m) be a local ring of dimension of n. Assume that the order
ideal conjecture is valid for local rings of dimension $(n-1)$. Then for the validity of the
order ideal conjecture on R it is enough to prove the validity of the assertion for the first
syzygies of modules of finite projection. In particular, for cyclic modules of finite projective
dimension over R, it is enough to prove that every minimal generator of any ideal in R of
finite projective dimension over R is a non-zero-divisor in R.
\endproclaim

\demo {Proof} Let $(F_\bullet, d_\bullet)$ be a minimal free resolution of M where $F_i = R^{r_i}$ for $i \ge 0$.
let $S_i$ denote the $i^{th}$ syzygy of M for $i \ge 1$ and let $\beta$ be a minimal generator of
$S_i$ for $i > 1$. Then $\beta = d_i(e)$ for some free generator $e$ of $F_i$ and we have
$\beta = \dsp \pmatrix a_1\\ \vdots \\ a_{r_{i-1}}\\ \endpmatrix \in \ R^{r_{i-1}}$.
Let $J$ denote the ideal generated by $a_1, \dots, a_{r_{i-1}}$.
Let $x$ be a non-zero-divisor on $R$; then $x$ is a non-zero-divisor on $S_1$.
Let $\overline R = R/xR$, $\overline S_i = S_i/xS_i$, $i \ge 1$. $\overline S_i$
is of finite projective dimension over $\overline R$ and
$\overline S_i = Syz^{i-1}_{\overline R} (\overline S_1)$ for $i > 1$. By induction hypothesis,
$\grade_{\overline R} (J+xR)/xR \ge (i-1)$; this implies $\grade_R J \ge (i-1)$.
Let $y \ \in J$ be a non-zero-divisor on R. Let $\tilde R = R/yR$,
$\tilde S_i = S_i/yS_i$ for $i \ge 1$. Then, again by arguing as above,
$\grade_{\tilde R} J/(y) \ge (i-1)$. Hence $\grade_R J \ge i$.
The second assertion now follows readily.
\enddemo

\proclaim{2.2\enspace Lemma} Let $M$ be a finitely generated module of finite projective dimension
over a local ring $(R, m)$. Suppose $\rank_R M = s$. Then there exists a free submodule $F = R^s$,
generated by a part of a minimal set of generators of $M$, such that $M/F$ has positive grade.
\endproclaim

\demo{Proof} We induct on $s = \rank_R M$. Since $\pd_R M < \infty$, if $s = 0$, then for every
associated prime $P$ of $R$, $M_P = 0$ and hence $\grade_R M > 0$. Now suppose that $s > 0$. Since
$\pd_R M < \infty$, by basic element method ( lemma 2.1, [E-G4]) there exists a minimal generator
$\alpha$ of $M$ such that image of $\alpha$ is a part of a basis of $M_P$ for every associated prime
$P$ of $R$. Hence we have a short exact sequence
$$
0 \to R \to M @>\phi>> M' \to 0  \ \ ( 1 \to \alpha)
$$
Then $\pd_R M' <\infty$ and $\rank_R M' = \rank_R M-1 = s-1$. By induction, there exists a free submodule
$F' = R^{s-1}$ of $M'$, generated by a part of a minimal set of generators of $M'$, such that a $M'' = M'/F'$
has positive grade. We have the following short exact sequence:
$$
0 \to F' @>i>> M' @>\psi>> M'' \to 0
$$
Since $F'$ is free we can lift $i : F' \to M'$ to $\eta : F' \to M$ such that $\phi . \eta = i$.
Let $\theta = \psi . \phi$. It can be easily checked that $\theta$ is onto and $\Ker \theta = R \oplus F'$.
Hence the lemma follows.
\enddemo

${\bold{2.3.}}$ Our next theorem reduces the order ideal conjecture to the assertion that every minimal
generator of a certain class of ideals of finite projective dimension must be a non-zero-divisor.

\proclaim{\enspace Theorem} Let $(R, m)$ be a local ring of dimension $n$. Assume that the order ideal
conjecture is valid for local rings of dimension $(n-1)$. Then for the validity of the order ideal
conjecture over $R$ it is enough to prove that every minimal generator of any ideal of grade 2,
height 2 and of finite projective dimension over $R$ is a non-zero-divisor in $R$.
\endproclaim

\demo{Proof} Let $M$ be a finitely generated module of finite projective dimension over $R$. Due to Lemma 2.1, for the validity of the order ideal conjecture it is enough to consider minimal generators of $S = Syz^1_R(M)$. If $\grade_RM = 0$, then $\ann_R M = 0$; Let $s_0 = \rank_R M$. By the previous lemma there exists a free submodule $F = R^{s_0}$ generated by a part of a minimal set of generators of $M$ such that $M' = M/F$ has positive grade. From the commutative diagram below
$$
\matrix & & & & 0 & & 0\\
& & & & \downarrow & & \downarrow\\
& & & & R^{s_0} & & R^{s_0}\\
& & & & \downarrow & & \downarrow\\
0 & \to & S & \to & R^{r_0} & @>\eta>> & M & \to & 0\\
& & \;|\wr & & \downarrow & & \downarrow\\
0 & \to & S & \to & R^{r_0-s_0} & \to & M^\prime & \to & 0\\
& & & & \downarrow & & \downarrow\\
& & & & 0 & & 0\endmatrix\tag{1}
$$
\enddemo

it follows that without any loss of generality we can assume $\grade_R M > 0$. Let $x$ be a non-zero-divisor in $\ann_R M$. By corollary to theorem 1.2 we can assume that $M$ has finite projective dimension over $R/xR$ and $\support(M) = \support (R/xR)$. Let $\overline R = R/xR$; $S = S/xS$. By tensoring $0 \to S \to R^{r_0} \to M \to 0$ with $\overline R$ we obtain the following short exact sequences
$$
0 \to M @>\theta>> \overline S \to T \to 0, \ \ 0 \to T \to {\overline R}^{r_0} \to M \to 0 \tag2
$$
where $T = Syz^1_{\overline R} (M)$. If $\gamma$ is a minimal generator of $M$ and $e \in R^{r_0}$, a free generator of $R^{r_0}$ is such that $\eta (e) = \gamma$, then it follows from chasing the commutative diagram obtained from
multiplying the short exact sequence $0 \to S \to R^{r_0} \to M \to 0$ by $x$ that
$\theta (\gamma) = \im(xe)$ in $\overline S$. Let {$\alpha_i \in S$} denote the lifts of a minimal set of generators {$\overline \alpha_i$}, $1\le i\le h$, of $T$. By induction, $\grade \ocal_{T, \overline R} (\overline{\alpha_i}) \ge 1$; hence $\grade \ocal_S (\alpha_i) \ge 1$ for $1 \le i \le h$. Due to the exact sequences in (2) we obtain a minimal set of generators $xe_1, \dots, x e_a, \alpha_1, \dots, \alpha_h$ of $S$ where $e_1, \dots e_a$ form a part of a basis of $R^{r_0}$. If $(F_\bullet, d_\bullet)$ $(F_i = R^{r_i})$ denote a minimal free resolution of $M$ over $R$, then there exist ${\tilde e}_1, \dots, {\tilde e}_a, \dots, {\tilde e}_{a+j}, \dots$ a basis of $R^{r_1}$ such that $d_1({\tilde e}_i) = x e_i$, $1 \le i \le a$ and $d_1 ({\tilde e}_{a+j}) = \alpha_j$, $1 \le j \le h$. Any minimal generator of $S$
is of the form $$\sum_{i=1}^t c_ixe_i + \sum_{j=1}^h d_j \alpha_j$$ where at least one of $c_i$s or $d_j$s is a unit. If any $d_j$ in the above expression is a unit then we are done by induction. Thus, it is easy to check that in order to show that for any minimal generator $\beta$ of $S$, $\grade \ocal_S (\beta) \ge 1$, it is enough to consider $\beta = x e - \sum{\lambda_i} \alpha_i$ where $e \in \{ e_1, \dots e_a \}$, $\sum \lambda_i \overline{\alpha}_i \ne 0$ in $T$. Due to the 2nd exact sequence in (2), if any $\lambda_i$ is a unit then we are done by induction. Hence we can assume that all $\lambda_i$s $\in m$ in the above expression of $\beta$.
Let $t = \rank_{\overline R} M$. By arguing as in (1) we obtain the following commutative diagram
$$
\matrix & & 0 & & 0 & & 0\\
& & \downarrow & & \downarrow & & \downarrow\\
0 & \to & R^{t} & @>x>> & R^{t} & \to & \overline R^{t}\\
& & \ \ \downarrow{\psi} & & \ \ \downarrow{\phi_0} & & \ \ \downarrow{\phi}\\
0 & \to & S & \to & R^{r_0} & \to & M & \to & 0\\
& & \downarrow & & \downarrow & & \downarrow\\
0 & \to & S^\prime & \to & R^{r_0-t} & \to & M^\prime & \to & 0\\
& & \downarrow & & \downarrow & & \downarrow\\
& & 0 & & 0 & & 0\endmatrix\tag{3}
$$

 Let $\{\overline{e}_j\}, 1 \le j \le t$ denote a basis of $\overline{R}^t$ such that $\gamma_j = \phi({\overline e}_j)$, $1 \le j \le t$, in (3) is a part of a minimal set of generators of $M$ (lemma 2.2). By commutativity of (3), $\psi (e_j) = x e_j$, $1 \le j \le t$. Note that none of $x e_j$ may be a part of $\{ x e_1, \dots, x e_a \}$ ---- part of a minimal basis of $S$ mentioned above. We want to prove the following:

\demo{Claim} For any $j$, $1 \le j \le t$, grade of the ideal generated by the entries of $x e_j - \sum \lambda_i \alpha_i$ is $\ge 1$ (here $xe_j = \psi (e_j))$.
\enddemo

\demo{Proof of the Claim} Let $\gamma_1, \dots, \gamma_t, \gamma_{t+1}, \dots, \gamma_{r_0}$ be a minimal set of generators of $M$ where $\gamma_j = \phi ({\overline e}_j)$, $1 \le j\le t$ as above, $t = \rank_{\overline R} M$. Then $\grade_{\overline R} ({M^\prime}) \ge 1$ i.e. $\grade_R ({M^\prime}) \ge 2$. Let $y \in \ann_{\overline R} {M^\prime}$ be a non-zero-divisor in ${\overline R}$; then $\forall k > t$, $y \gamma_k = \sum^t_{j=1} a_{kj} \gamma_j$; hence $ye_k - \sum a_{kj} e_j \in S$ for $k > t$. Let $P \in Ass_R (R)$ and let $q \in Ass_{\overline R} ({\overline R})$ containing $P$. By construction $q$ is an associated prime of $M$ (recall $\ann_{\overline R} M = 0$); then $M_q$, $T_q$ are free ${\overline R}_q$ modules of $\rank t$ and $r_0 - t$ respectively. Since $\pd_{R_q} M_q = 1$, $S_q$ is also a free
$R_q$-module. We have the following short exact sequence
$$
0 \to M_q \to {\overline S}_q \to T_q \to 0   \ \ \ \ \ \  0 \to T_q \to {\overline R}_q^{r_0} \to M_q \to 0 \tag4
$$
Let ${\overline \beta}_k = {\overline e}_k - \sum_{j=1}^t \dfrac {{\overline a}_{kj}} {\overline y} {\overline e}_j$. Then $\{ {\overline \beta}_k \}, {t+1} \leq k \leq r_0$ is a basis of $T_q$. Let $\beta_k = e_k - \sum_{j=1}^k \dfrac {a_{kj}} {y} e_j, t+1 \leq k \leq r_0$. Then it follows from (4) $\{xe_1 \dots, xe_t, \beta_{t+1}, \dots, \beta_{r_0}\}$ form a basis of $S_q$. In $S_q$, we have
$$
\sum \lambda_i \alpha_i = \sum_{i=1}^t \dfrac {c_i} {b} xe_i + \sum_{k=t+1}^{r_0} \dfrac {d_k} {b} \beta_k,  b  \not  \in q \tag5
$$
If $\sum \lambda_i \alpha_i \not \in qS_q$, then $\sum \lambda_i {\overline \alpha_i} \not \in qT_q$ which implies that some $d_k \not \in q$. Hence $\im (xe_j - \sum \lambda_i \alpha_i)$ is a minimal generator of $T_q$ which is free and thus $xe_j - \sum \lambda_i \alpha_i \not \in PR^{r_0}$. Now suppose that $\sum \lambda_i \alpha_i \in qS_q$ then all $c_i$s, $d_i$s $\in q$ and $x \nmid d_k$ for at least one $k$ in (5). Hence $xe_j - \sum \lambda_i \alpha_i = (1- \dfrac {c_j} {b}) xe_j - \sum_{i \neq j} \dfrac {c_i} {b} xe_i - \sum \dfrac {d_{ik}} {b} \beta_k$. Since $1-\dfrac {q} {b}$ is a unit in $R_q$, $xe_j - \sum \lambda_i \alpha_i$ is a minimal generator of $S_q$. Since $S_q$ is a free $R_q$ module, $xe_j - \sum \lambda_i \alpha_i \not \in PR_q^{r_0}$ and hence $xe_j - \sum \lambda_i \alpha_i \not \in PR^{r_0}$. This  completes the proof of our claim.

Due to the commutative diagram (3) and the above Lemma we can assume without any loss of generality $\grade_R M \ge 2$. Now we appeal to the following result due to Smoke (Lemma 4.1, Th. 4.2, [Sm]). Given a finitely generated $R$-module $M$ of $\grade \ge 2$, we can construct an exact sequence
$$
0 \to M \to R/I \to R/J \to 0 \tag6
$$
where $R/J$ has a filtration whose successive quotients are isomorphic to cyclic modules of the form $R/(u, v)$ where $\{u, v\}$ form an $R$-sequence. The proof in [Sm] shows that, by choosing the $R$-sequences of length 2 in $m$ (annihilator of the corresponding module over $R$), this exact sequence can be constructed in such a way that if $\delta:S \to I$ denotes the corresponding map on first syzygies via (6), then $\delta$ is surjective and ${\overline \delta} : S/mS \to I/mI$ is an isomorphism. Since $\pd_RM < \propto$, $\pd_R R/I < \propto$ and it follows from the construction that height of $I = $ grade of $I = 2$. Thus the proof of our theorem is complete.
\enddemo

${\bold{2.4.}}$ Now we state the following Lemma.
\proclaim {Lemma} Let (R, m) be a local ring and let $M$ be a finitely generated
module of projective dimension $n < \infty$. Let $N$ be an $R$-module such that
$\ann_R M$ contains an $N$-sequence of length $r$. Then $\Tor_{n-i}^R (M, N) =  0$ for
$0 \le i < r$.
\endproclaim

We leave the proof of this lemma to the reader.

${\bold{2.5.}}$ In our next theorem we indicate the vanishing of $\Tor$ from an altogether
different perspective.

\proclaim{Theorem}Let $(R,m)$ be a local ring and let
$M$ be a finitely generated module of finite projective dimension over
$R$. Let $\pd_RM=h$. Then there exists an $R$-sequence $x_1,\dots,x_h$
of length $h$ such that for any $R$-module $N$ for which $x_1,\dots,x_h$
form an $N$-sequence, $\Tor_i^R(M,N)=0$ for $i>0$.
\endproclaim

\demo{Proof}If $\rank_RM=s_0>0$, by lemma 2.2, there
exists a free submodule $F=R^{s_0}$ generated by a part of a
minimal set of generators of $M$ such that $M' = M/F$
has positive grade. From the commutative diagram
$$
\matrix & & & & 0 & & 0\\
& & & & \downarrow & & \downarrow\\
& & & & R^{s_0} & & R^{s_0}\\
& & & & \downarrow & & \downarrow\\
0 & \to & S & \to & R^{r_0} & \to & M & \to & 0\\
& & \;|\wr & & \downarrow & & \downarrow\\
0 & \to & S & \to & R^{r_0-s_0} & \to & M^\prime & \to & 0\\
& & & & \downarrow & & \downarrow\\
& & & & 0 & & 0\endmatrix\tag{1}
$$
it follows that $S=\Syz^1_R(M) = \Syz_R^1(M^\prime)$. If $\rank_R M = 0$, then $M' = M$.

Let $x_1\in \ann_RM^\prime$ be a non-zero-divisor on $R$ and let $R_1=R/xR$,
$S_1=S/xS$. By theorem 1.2 we have an exact sequence of $R_1$-modules
$$
0\to M^\prime \to M_1 \to V_1 \to 0\tag{2}
$$
such that $\pd_{R_1}M_1<\infty$ $(=\pd_RM-1)$ and
$\pd_{R}V_1=1$. It is also clear from Corollary 1.3 that
$S(M_1)=\Syz^1_R(M_1)=S\oplus R^{t_1} = \Syz^1_R (M') \oplus R^{t_1} = \Syz^1_R (M) \oplus R^{t_1}$
for some $t_1\ge 0$. Tensoring the short exact sequence $0\to S(M_1)\to R^{\ell_1}\to
M_1\to 0$ with $R_1$, we obtain the following short exact sequence
$$
0\to M_1\to S(M_1)\otimes R_1\to T_1\to 0,\quad 0\to T_1\to
R_1^{\ell_1}\to M_1\to 0.\tag{3}
$$
Here $T_1=\Syz^1_{R_1}(M_1)$. Recall that $\pd_{R_1}M_1<\infty$ and
$\support(M_1)=\spec R_1$. Now we start with $M_1$ over $R_1$ and
repeat the process described in (1), (2) and (3). We continue this
process $(h-2)$ times and obtain an $R$-sequence $x_1,\dots,x_i$, $1\le
i\le h-1$, modules $M_i$ of finite projective dimension over
$R_i=R/(x_1,\dots,x_i)$ and short exact sequences
$$
\align
&0\to M_{i-1}^\prime\to M_i\to V_i\to 0\tag{2i}\\
&0\to F_i\to M_i\to M_i^\prime\to 0\tag{1i}
\endalign
$$
where $\pd_{R_i}M_i=\pd_RM-i$, support$(M_i)$ = support$(R_i)$, $\pd_{R_{i-1}}V_i=1$,
$M_{i-1}^\prime = M_{i-1}/F_{i-1}$, $F_{i-1}$ a free $R_{i-1}$-module
as constructed in (1), $M_i^\prime$ is a module over $R_{i+1}$. We
also have short exact sequences of $R_i$-modules
$$
0 \to M_i\to S(M_i)\otimes R_i\to T_i\to 0,\;\;0\to T_i\to
R_i^{\ell_i}\to M_i\to 0.\tag{3i}
$$
We note that grade $M_i=i$ and $\pd_RM_i=h$. Let $x_h\in m$ be such
that $im(x_h)$ in $R_{h-1}$ is a non-zero-divisor contained in
$\ann_{R_{h-1}}M^\prime_{h-1}$. Then $x_1,\dots,x_h$ form an $R$-sequence
and this is our required sequence.

Since $\pd_R M'_{h-1} = h$, if $N$ is an $R$-module such that $x_1,\dots,x_h$ form an
$N$-sequence, then it follows from the above lemma that $\Tor_i^R (M_{h-1}^\prime, N) =  0$
for $i > 0$. Now it follows from the above short exact sequences
starting from $M^\prime_{h-1}$ and tracing back to $M$ that
$\Tor^R_j(M,N)=0$ for $j\ge 1$. If $\grade_RM=r>0$ and $x_1,\dots,x_r$
is an $R$-sequence contained in $\ann_RM$, we start with $Q$---an
$\overline{R}(=R/(x_1, \dots, x_r))$ module as in Theorem 1.2 and construct an
$\overline{R}$-sequence $x_{r+1},\dots,x_h$ by the above method. Then
$x_1,\dots,x_h$ form an $R$-sequence satisfying the required vanishing
property of $\Tor$.

Since $\pd_RM < \infty$, if $\grade_RM = 0$ then $\ann_RM =0$ and hence $\rank_RM > 0$.
Thus we are back to diagram (1).
\enddemo

\proclaim{2.6 \enspace Corollary 1}Let $(R,m)$ be a local ring and let $M$ be a
finitely generated module of finite projective dimension $h$ over $R$. Let
$x_1,\dots,x_h$ be an $R$-sequence as mentioned in the above
proposition. Suppose that for every $P\in \Ass_R(R)$, there exists an
$R/P$-module $N$ such that $x_1,\dots,x_h$ form an $N$-sequence and
$N\ne mN$. Then, for every minimal generator $\beta$ of
$S_1=\Syz^1_R(M)$, $\grade\ocal_{S_1}(\beta)\ge 1$.
\endproclaim

\demo{Proof}Let $(F_\bullet,d_\bullet)$ be a minimal resolution of $M$
over $R$. If possible let grade $\ocal_{S_1}(\beta)=0$. Then there
exists an associated prime $P$ such that  $\ocal_{S_1}(\beta)\subset
P$. Let $\overline{R}=R/P$ and $\overline{F}_\bullet =
F_\bullet\otimes R/P$. Consider the sequence
$$
\overline{F_2} @>\overline{d_2}>> \overline{F_1} @>\overline{d_1}>>
\overline{F_0}\to 0.
$$
Since $\ocal_{S_1}(\beta)\subset P$, $\overline{d}_1(\beta)=0$. This
implies that for any $\overline{R}$-module $N$, $\Tor_1^R(M,N)\ne
0$. However, by hypothesis, there exists an $\overline{R}$-module $N$
such that $x_1,\dots,x_h$ form an $N$-sequence. Then, by the above
proposition, we have $\Tor^R_j(M,N)=0$ for $j>0$, which leads to a
contradiction. Hence $\grade\ocal_{S_1}(\beta)\ge 1$.
\enddemo

\proclaim{2.7 \enspace Corollary 2}Let $(R,m)$ be a local ring and $M$ be a
finitely generated module of finite projective dimension $h$. Let
$x_1,\dots,x_h$ be an $R$-sequence as mentioned in the above
proposition. If $R$ has mixed characteristic $p>0$, we assume that
$p,\,x_1,\dots,x_h$ form a part of a system of parameters of $R$. Then
for every minimal generator $\beta$ of $S_1=\Syz_R^1(M)$, $\grade
\ocal_{S_1}(\beta)\ge 1$.
\endproclaim

\demo{Proof}If possible, let $\grade\ocal_{S_1}(\beta)=0$. Then there
exists an associated prime $P$ of $R$ such that
$\ocal_{S_1}(\beta)\subset P$. If $R$ is equicharacteristic then
there exists a big Cohen-Macaulay $R/P$-module $N$ such that
$x_1,\dots,x_h$ form a regular $N$-sequence. If $R$ has mixed
characteristic $p>0$, then there exists a big Cohen-Macaulay
$R/(P+pR)$-module $N$ such that $x_1,\dots,x_h$ form a
regular $N$-sequence. Hence we are done by Corollary~1.
\enddemo

${\bold{Remark}}$ Due to the existence of Cohen-Macaulay algebras over
local domains of dimension less than or equal to three ([H3]) it can be checked
from the above arguments that finitely generated modules of projective dimension
less than or equal to three satisfy the order ideal conjecture.

\proclaim{2.8 \enspace Lemma}Let $(R,m,K=R/m)$ be an
equicharacteristic complete local ring of dimension $d$ and let
$x_1,\dots,x_d$ be a system of parameters of $R$. Let $M$ be a big
Cohen-Macaulay $R$-module such that $\bx M\ne M$ and $x_1,\dots,x_d$
form a maximal $M$-sequence. Then $\widehat{M}\;=$ the $m$-adic
completion of $M$ is a flat $K[[x_1,\dots,x_d]]$-module.
\endproclaim

\demo{Proof}Let $S=K\[[x_1,\dots,x_d]\]$. Then $S$ is a complete
power series ring in $d$ variables, $R$ is a module finite extension
of $S$ and $\bx=\{x_1,\dots,x_d\}$ form a regular system of parameters
of $S$. Moreover $\bx$ is $\widehat{M}$-regular (Th.\ 8.5.1
\cite{B-H}) and $\widehat{M}$ is a balanced big Cohen-Macaulay module
(Cor.\ 8.5.3, \cite{Br-H}). Hence
$\dsp\psi:\dfrac{\widehat{M}}{\bx \widehat{M}} \[X_1,\dots,X_d\] \to
\bigoplus\limits^\infty_{n=0}
\dfrac{\bx^n\widehat{M}}{\bx^{n+1}\widehat{M}}$ is an
isomorphism. Since $\bx$ is a regular system of parameters of $S$,
$\dsp K\[X_1,\dots,X_d\] \simeq \bigoplus\limits^\infty_{n=0}
\dfrac{\bx^nS}{\bx^{n+1}S}$. Hence $\dsp\dfrac{\bx^nS}{\bx^{n+1}S}
\otimes_K \widehat{M}/\bx \widehat{M} \simeq
\dfrac{\bx^n\widehat{M}}{\bx^{n+1}\widehat{M}}$ and
$\dfrac{\widehat{M}}{\bx\widehat{M}}$ is a non-null vector spacing
over $K = S/\bx S$. Thus it follows, by Th.~1, \S 5.2 in \cite{Bou},
that $\widehat{M}$ is $S$-flat.
\enddemo

\proclaim {2.9 \enspace Theorem }{\rm (Foxby, \cite{F})}Let $(R,m,K)$ be
an equicharacteristic complete local ring and let $N$ be a finitely
generated module of finite projective dimension over $R$. Let $I$ be
an ideal of height $h$ and let $M$ be a big Cohen-Macaulay module over
$R/I$. Let $\widehat{M}$ be the $m$-adic completion of $M$. Then
$\Tor_i^R(\widehat{M},\,N)=0$ for $i>h$.
\endproclaim

\demo{Proof}For a proof we refer the reader to (\cite{E-G4}) or
(\cite{F}) where the existence of a big Cohen-Macaulay $R/I$-module
$Q$ such that $\Tor_i^R(Q,N)=0$, for $i>h$ has been demonstrated (it
was first proved by Foxby). In these proofs, it was required that such
a $Q$ be free over a certain complete regular local ring $S$ contained
in $R/I$ such that $R/I$ is a module-finite extension of $S$. By
Lemma 2.8 the completion $\widehat{M}$ of any big Cohen-Macaulay
$R/I$ module is flat over certain complete regular subrings of
$R/I$. And this flatness is enough to ensure the validity of arguments
provided in theorem 1.11 in (\cite{E-G4}) or in (\cite{F}) for proving our assertion.
\enddemo


\subhead{2.10}\endsubhead Next we recall Shimomoto's theorem.

\proclaim{Theorem {\rm (Th.\ 5.3, \cite{Shi})}}Let $(R,m)$ be a complete
local domain of mixed characteristic $p>0$. Then there exist some
system of parameters $p,x_2,\dots,x_d$ of $R$ and an almost
Cohen-Macaulay quasi-local algebra $B$ over $R^+$, the integral closer of $R$ in
algebraic closer of the field of fractions of $R$, in the sense that

\roster
\item"1." $(p,x_2,\dots,x_d)B\ne B$,

\item"2." $x_2,\dots,x_d$ form a regular sequence on $B/pB$, and

\item"3." $p$ is not nilpotent in $B$ and the ideal $(0:p)_{\bold{B}}$
is annihilated by $p^\in$ for any rational $\in > 0$.
\endroster
\endproclaim

Actually Shimomoto's construction shows that such a $B$ can be constructed for any system of
parameters of the form $p,x_2,\dots,x_d$ of $R$.

\smallskip

\noindent {\bf Definition.}\enspace An almost Cohen-Macaulay algebra $B$ as
above is called balanced if $B/pB$ is a balanced Cohen-Macaulay algebra.

\smallskip

\proclaim{2.11} \rm Now we are ready to prove our final theorem. \endproclaim

\proclaim{\enspace Theorem}Let $(R,m)$ be a local ring of mixed
characteristic $p>0$. We assume that either $p$ is nilpotent or $p$ is
a non-zero-divisor in $R$. Let $M$ be a finitely generated module of
finite projective dimension over $R$ and $\beta$ be a minimal
generator of $S_i$, the $i^{\text{th}}$ syzygy of $M$ (minimal), for
$i>0$. We have the following:

\roster
\item"a)" If $p$ is nilpotent, the order ideal conjecture is valid on~$R$;

\item"b)" if $pM=0$, then grade of $\ocal_{S_1}(\beta)\ge 1$; and

\item"c)" Suppose that balanced or complete almost Cohen-Macaulay algebras exist
over complete local domains. If $pM=0$ and $p\in m-m^2$, then grade of
$\ocal_{S_i}(\beta)\ge i$, $\forall i\ge1$.

\item"d)" Assume that every element in $m-m^2$ is a non-zero-divisor and that the order
ideal conjecture is valid over  $R/xR$ for any $x \ \in \ m-m^2$. If $\ann_RM \cap (m-m^2) \ne \phi$
or $depth_R M > 0$ then $grade\ocal_{S_i}(\beta)\ge i$, $\forall i \ge 1$.

\endroster
\endproclaim

\demo{Proof}\newline

a)\enspace If $p$ is nilpotent, the proof follows immediately by
similar arguments as in (Th. 2.4, \cite{E-G2}) due to the existence of big
Cohen-Macaulay modules on equicharacteristic local domains $R/P$ for
every prime ideal $P$ in $\spec R$. One may also use Lemma 9.1.8 from
\cite{Br-H} to prove the assertion.

b)\enspace Since $p$ is a non-zero-divisor on $R$, by the corollary of Theorem
1.2 we can assume $\pd_{R/pR}M<\infty$. We write
$\overline{R}=R/pR$. Consider the short exact sequence
$$
0\to S_1 @>g>> R^{r_0} @>\eta>> M\to 0.
\tag{1}
$$
Tensoring this sequence with $\overline{R}$ we obtain the following
exact sequences
$$
\align
&\ocal\to M @>j>> \overline{S}_1 @>b>> T_1\to 0,\tag{2}\\
&\ocal\to T_1 @>\gamma>> \overline{R}^{r_0} @>\overline{\eta}>>
M\to 0\tag{3}
\endalign
$$
where $T_1=\Syz^1_{\overline{R}}(M)$, $\overline{S}_1=S_1\otimes_R
\overline{R}$. For any minimal generator $\vartheta$ of $M$ we have
$\vartheta=\eta(e)$, $e$ being a minimal generator of $R^{r_0}$
and $j(\vartheta)=$ class of $pe$ in $\overline{S}$. Recall that if
$I$ is an ideal of $R$ such that $\grade_RI=0$, then for any
non-zero-divisor $x$ in $R$, $\grade_{\overline{R}}(I+xR)/xR$ is also
$0$. Hence, for any minimal generator $\beta$ of $S_1$, if
$b(\overline{\beta})$ is a minimal
generator of $T_1$ then our assertion follows due to the validity of
the order ideal conjecture on equicharacteristic local rings. If
$b(\overline{\beta})=0$, then $\overline{\beta}\in j(M)$ i.e.\
$\overline{\beta}=j(\vartheta)$ where $\vartheta$ is a minimal generator of
$M$. Let $\{\overline{\alpha}_j\}_{1\le j\le t}$ be a minimal set of
generators of $T_1$ and $\{\alpha_j\}$ denote a lift of
$\{\overline{\alpha}_j\}$ in $S_1$. Then there exists a basis
$e_1,\dots, e_i\dots, e_{r_0}$ of $R^{r_0}$ such that
$\{p\,e_1,\dots,p\,e_i,\;\alpha_js\}$ form a minimal set of
generators of $S_1$. Moreover, in order to show that for any minimal
generator $\beta$ of $S_1$, $\grade\ocal_{S_1}(\beta)\ge 1$, it is
enough to take $\beta=p\,e-\sum\limits^t_{j=1} \lambda_j\alpha_j$, where
$e\in \{e_1,\dots,e_i\}$, $\lambda_js\in m$ and
$\sum{}\lambda_j\overline{\alpha}_j\ne 0$ in $T_1$.

If possible let $\grade\ocal_{S_1}(\beta)=0$ i.e.\ $\ocal_{S_1}(\beta)\subset
P$, for some $P\in \Ass_R(R)$. Let $\pd_RM=h$; then
$\pd_{\overline{R}}M=(h-1)$. By theorem 2.5, corresponding to $M$
over $\overline{R}$, there exists an $\overline{R}$ sequence
$\overline{x}_1,\dots,\overline{x}_{h-1}$ satisfying the assertion mentioned in
theorem (2.5). Let $\widetilde{R}=R/PR$
and let $\tilde{p}=im(p)$ in $\widetilde{R}$ and
$\tilde{x}_i=im\,(x_i)$ in $\widetilde{R}$, $1\le i\le h-1$. Then
$\tilde{p}$, $\tilde{x}_1,\dots,\tilde{x}_{h-1}$ form a part of a system of
parameters of $\widetilde{R}$. By Corollary 2.7 there exists an
almost Cohen-Macaulay $\widetilde{R}^+$ algebra $B$ such that
$(\tilde{p},\dots,\tilde{x}_{h-1})B\ne B$ and $\tilde{x}_1,
\dots,\tilde{x}_{h-1}$ form a regular sequence on $B/pB$. Then, by
theorem 2.5 or by theorem 2.4 in [F], we have $\Tor_j^{\overline{R}}(M,\,B/pB)=0$ for
$j>0$.

We consider the following part of a minimal resolution $F_\bullet$ of
$M$ over $R$:
$$
\matrix
R^{r_2} & @>d_2>> & R^{r_1} & @>d_1>> & R^{r_0} & \to & 0\\
\| & & \| & & \|\\
F_2 & \to & F_1 & \to & F_0 & \to & 0
\endmatrix.
$$
Let $\widetilde{F}_i = F_i\otimes \widetilde{R}$, $0\le i\le
2$. Tensoring the above sequence with $\widetilde{R}$ we obtain a
sequence
$$
\widetilde{F}_2 @>\tilde{d}_2>> \widetilde{F}_1 @>\tilde{d}_2>>
\widetilde{F}_0 \to 0
$$
where $\tilde{\beta}=im(\beta)$ in $\widetilde{F}_0 = 0$.

Tensoring the above sequence with $B$ and writing $F_{iB}=F_i\otimes B$ we have
$$
F_{2B} @>d_{2B}>> F_{1B} @>d_{1B}>> F_{0B}\to 0\tag{4}
$$
where $\beta_B=im\tilde{\beta}$ in $F_{0B}=0$. Hence
$$
pe_B=\sum{}\lambda_j\alpha_{jB}\quad\text{in}\quad
B^{r_0}=F_{0B},\tag{5}
$$
where $\alpha_{jB}=im(\alpha_j\otimes 1_B)$ in $F_{0B}$, $\alpha_j\otimes
1_B\in S_1\otimes B$. Since $\Tor_j^{\overline{R}}(M,B/pB)=0$, for
$j>0$, tensoring (2) and (3) with $B/pB$ over $\overline{R}$ we
obtain the following short exact sequences:
$$
0 \to M\otimes B/pB @>j \bigotimes 1_{B/pB}>> \overline{S_1}\otimes B/pB\to T_1\otimes B/pB\to 0\tag{6}
$$
and
$$
\ocal\to T_1\otimes B/pB\to (B/pB)^{r_0}\to M\otimes B/pB\to 0.\tag{7}
$$
Let $\overline{\alpha_{jB}}$ $=$ the image of
($\overline{\alpha_j}\otimes 1_{B/pB}$ in ${T_1}\otimes B/pB$) in
$\(B/pB\)^{r_0} = im\,\alpha_{jB}$ in $(B/pB)^{r_0}$. Due to (5) and
(7), we have $\sum{}\lambda_j\overline{\alpha_{jB}}=0$ in
$(B/pB)^{r_0}$ and hence $\sum{}\lambda_j(\overline{\alpha_j}\otimes
1_{\overline{B}})=0$ in ${T_1}\otimes B/pB$. Since $\lambda_js \in\/m$, this
implies, due to the exact sequence~(6) and the definition of $j$ in ~(2), that
$$
\sum{}\lambda_j(\alpha_j\otimes 1_B) = \sum{} a_i(p\,e_i\otimes 1_B) +
p\(\sum{}b_i(p\,e_i\otimes 1_B) + \sum{}\mu_i(\alpha_i\otimes 1_B)\)\tag{8}
$$
in ${S_1}\otimes B$, $a_is\in m_B$, where $m_B =$ maximal ideal of $B$.

Hence, we have from (5) and (8).
$$
p\,e_B=\sum{}(a_i+pb_i)\,p\,e_{iB}+p\sum{}\mu_i\alpha_{iB}
$$
in $B^{r_0}$ i.e.
$$
p\left\{e_B-\[\sum{}(a_i+p\,b_i)e_{iB} +
\sum{}\mu_i\alpha_{i\beta}\]\right\}=0.\tag{$\ast$}
$$
Since $a_i$, entries of $\alpha_{i\beta}\in m_B$,
$e_B-\[\sum{}(a_i+p\,b_i)e_i + \sum{}\mu_i\alpha_{iB}\]$ is a free
generator of $B^{r_0}$ and hence $p$ can not annihilate it. Thus $(\ast)$
leads to a contradiction. Hence $\grade\ocal_{S_1}(\beta)$ must be $\ge 1$.

\smallskip

c)\enspace We assume $p\in m-m^2$ and $pM=0$. Let $F_\bullet : 0\to
R^{r_n} @>d_n>>
R^{r_{n-1}} \to \cdots \to R^{r_1} @>d_1>> R^{r_0} \to
0$ be a minimal free resolution of $M$ over $R$ and let $P_\bullet :
0\to \overline{R}^{s_{n-1}}\to \cdots \to \overline{R}^{s_1} \to
\overline{R}^{s_0} \to 0$ be a minimal free resolution of $M$ over
$\overline{R}(=R/pR)$. Shamash (\cite{Sha}) has shown that $P_\bullet$
can be obtained from $F_\bullet$ via the homotopy maps $h_\bullet :
F_\bullet \to F_\bullet (+1)$ induced by the $0$-map on $M$ due to
multiplication by $p$. Actually each $F_i$, for $i>0$, can be
decomposed into two parts: $F_i=F_i^\prime\oplus F_i^{\prime\prime}$,
where $h_{i}(F_{i}^{\prime\prime})=0$, every free generator
$e^\prime$ of $F_{i-1}^\prime$ is such that $e=h_{i-1}(e^\prime)$ is a
free generator of $F_i^{\prime\prime}$ and
$d_i(e)=p\,e^\prime-h_{i-2}d_{i-1}(e^\prime)$. Moreover, it follows from Shamash's
theorem that I) there exists $\phi_\bullet : P_\bullet\to
\overline{F_\bullet(+1)}(\overline{F}_\bullet = F_\bullet\otimes_R
\overline{R})$, where $\phi_0$ induces the inclusion map:
$M=\Tor_1^R(M,\,R/p){\;\overset{j}\to{\hookrightarrow}\;}\overline{S}_1$ and
II) $\phi_\bullet$ induces a splitting on each
component of $P_\bullet$ and $P_\bullet(+1)$ can be extracted from the
mapping cone of $\phi_\bullet$. Let $S_i=\Syz^i_R(M)$, $T_i=\Syz^i_{\overline{R}}(M)$ and
$\overline{S_i}=S_i\otimes \overline{R}$. This leads to the following short exact
sequences for $i>1$:
$$
\align
&0\to T_{i-1} @>\psi>> \overline{S}_i @>\eta>> T_i\to 0,\tag{$1^\prime$}\\
&0\to T_i\to \overline{R}^{r_{i-1}}\to T_{i-1}\to 0\tag{$2^\prime$}
\endalign
$$
where $\psi$ is induced by $\phi_\bullet$ and
$\psi(\im\,\overline{e}^\prime) =
\overline{p\,e-h_{i-2}\,d_{i-1}(e^\prime)}$ is a minimal generator of
$\overline{S_i}$ due to the splitting property of $\phi_\bullet$. Let
$\gamma=p\,e^\prime - h_{i-2}d_{i-1}(e^\prime)$ and
$$
\overline{\gamma} = \overline{p\,e^\prime - h_{i-2}\,d_{i-1}(e^\prime)}.\tag{$3^\prime$}
$$
Then $\gamma,\overline{\gamma}$ are minimal generators of $S_i$ and
$\overline{S_i}$ respectively.
\enddemo 

\medskip

\noindent {\bf Claim:}\enspace $\grade\ocal_{S_i}(\gamma)\ge i$.

\medskip

\demo{Proof of the claim}If possible let $P$ be a prime ideal of
$\height (i-1)$ containing $\ocal_{S_i}(\gamma)$. Since $e^\prime\in F_{i-1}^\prime$
and $h_{i-2}(F_{i-2}) = F_{i-1}^{\prime\prime}$, it follows that $p\in P$. Let
$p,\,x_2,\dots,x_{i-1}$ be a maximal $R$-sequence contained in
$\ocal_{S_i}(\gamma)$---we denote it by $\bx$. Then
$\Tor_j^R(M,\,R/\bx)=0$ for $j\ge i$. Let $R^\prime = R/\bx$,
$S_i^\prime = S_i\otimes R^\prime$ etc. We have an exact sequence
$$
0\to S_i^\prime \to R^{\prime\,r_{i-1}} \to S^\prime_{i-1} \to 0.
$$
Then $S^\prime_{i-1}$ is a module of finite projective dimension over
$R^\prime$ and $S_i^\prime=\Syz^1_{R^\prime}\,(S^\prime_{i-1})$
has a minimal generator $\gamma^\prime = im\,\gamma$ in $S_i^\prime$ such
that $\ocal_{S_i^\prime}(\gamma^\prime)$ has grade $0$ in $R^\prime$.
This contradicts part b) of our theorem and hence the claim is
established.

Let $\{\overline{\alpha}_j\}$ form a minimal set of generators of
$T_i$ and let $\{\alpha_j\}$ denote their lifts in $S_i$. Since
characteristic of $\overline{R}=p>0$, by Evans-Griffith theorem
(\cite{E-G2}),
$\grade_{\overline{R}}(\ocal_{T_i}(\overline{\alpha}_j))\ge i$. Since
$(1^\prime)\otimes R/m$ is exact, due to the above claim, for the purpose of
proving our theorem it would be enough to establish that for every
minimal generator $\beta$ of $S_i$ of the form
$\beta=\gamma-\sum{}\lambda_j\alpha_j$,
$\grade\ocal_{S_i}(\beta)\ge i$, where $\overline\gamma = \psi(im\, \overline{e^\prime})$ for some
free generator ${e^\prime}$ for $F_{i-1}$ and $\lambda_js\in {m}$. If possible let
$\ocal_{S_i}(\beta)\subset P-$ a prime ideal of height $(i-1)$. If
$p\in P$, there exists a maximal Cohen-Macaulay $R/P$-algebra
($(R/P)^+$ algebra) $\overline{B}$. If $p\not\in P$, then, by assumption, there exists
a balanced almost Cohen-Macaulay $R/P$-algebra ($(R/P)^+$ algebra) $B$
such that $p$ is not nilpotent on $B$ and $\overline{B}=B/pB$ is
maximal Cohen-Macaulay algebra over $R/(P+pR)$. Hence, in either case,
by Theorem 2.4 in [F], $\Tor_j^{\overline{R}} (M,\overline{B})=0$ for $j\ge
i$. Tensoring ($1^\prime$) and ($2^\prime$) with $\overline{B}$
we get the following exact sequences
$$
\ocal\to T_{i-1}\otimes \overline{B} @>\psi\otimes 1_{\overline{B}}>> \overline{S}_i\otimes
\overline{B}\to T_i\otimes \overline{B}\to 0\tag{$6^\prime$}
$$
and
$$
\ocal\to T_i\otimes \overline{B}\to \overline{B}^{r_{i-1}}\to
T_{i-1}\otimes \overline{B}\to 0.\tag{$7^\prime$}
$$
Let $e$ be a free generator of $F_i$ such that $d_i(e)=\beta$.

From $F_\bullet$, we consider the exact sequence
$$
F_{i+1} @>d_{i+1}>> F_i @>d_i>> F_{i-1}.
$$
Let $\widetilde{R}=R/P$, $\widetilde{F}_i=F_i\otimes_R
\widetilde{R}$. Tensoring the above sequence with $\widetilde{R}$, we
obtain a complex
$$
\widetilde{F}_{i+1} @>\tilde{d}_{i+1}>> \widetilde{F}_i
@>\tilde{d}_i>> \widetilde{F}_{i-1}
$$
where $\tilde{d}_i\;(\tilde{e})=0$. Let $e_B=e\otimes 1_B$. Then in
the complex $F_{i+1}\otimes B\to F_i\otimes B\to F_{i-1}\otimes B$, we
have $d_{iB}(e_B)=0$. Let $\gamma_B=im\,(\gamma\otimes 1_B)$,
$\alpha_{jB}=im\,(\alpha_j\otimes 1_B)$ in $B^{r_{i-1}}=F_{i-1}\otimes
B$. Then
$$
d_{iB}(e_B)=0 \Rightarrow
\gamma_B-\sum{}\lambda_j\,\alpha_{jB}=0.\tag{$5^\prime$}
$$
Let $\overline{\alpha}_{jB}=\im \alpha_{jB}$ in
$\overline{B}^{r_{i-1}}$. Due to $(5^\prime)$ and $(6^\prime)$ we have
$\sum{}\lambda_j\overline{\alpha}_{jB}=0$ in
$\overline{B}^{\gamma_{i-1}}$, hence
$\sum{}\lambda_j(\overline{\alpha}_j\otimes 1_B)=0$ in $T_i\otimes
\overline{B}$. Since $\lambda_js\in\/m$, this implies, due to the exact
sequence $(6^\prime)$ and definition of $\psi$ in $(1^\prime)$, that
$$
\sum{}\lambda_j(\alpha_j\otimes 1_B) = \sum{}a_i(\gamma_i\otimes 1_B)
+ p\[\(\sum{}\,b_i(\gamma_i\otimes 1_B) +
\sum{}\,\mu_j(\alpha_j\otimes 1_B)\)\]\tag{$8^\prime$}
$$

in $S_i \bigotimes\/ B, a_is \in m_B$. Hence from $(5^\prime)$ we have

$$
\gamma_B = \sum{}(a_i+pb_i)\;(\gamma_{iB}) + p\sum{}\,\mu_j\alpha_{jB}
$$
in $B^{r_i-1}$.
Since $\gamma_B=p\,e^\prime_B - h_{i-2}\,d_{i-1}(e^\prime_B)$ and
$\gamma_{iB}s$ also have similar expressions, comparing the
$e_B^\prime${}${\text{th}}$ co-ordinate in $B^{r_{i-1}}$, we obtain from above
$$
p\[e^\prime_B - \sum{}(a_i + pb_i)\,e_{B^\prime} - \sum{}\delta_j\] = 0\tag{**}
$$
where $\delta_j = e_B^\prime${}$^{\text{th}}$ co-ordinate of
$\mu_j\alpha_{jB}$. Since $a_i,\,\delta_j\in m_B$, the term within
brackets in (**) is a free generator of $B^{r_{i-1}}$ and hence $p$
cannot annihilate it. Thus (**) leads to a contradiction. Hence
$\grade\ocal_{S_i}(\beta)\ge i$ and our proof is complete.

d) We assume that every element in $m-m^2$
is a non-zero-divisor and for any such element x the order ideal conjecture is valid for $R/xR$.
Let M be a finitely generated module of finite
projective dimension such that either $\ann_R M \cap (m-m^2) \ne \phi$ or $depth_R M > 0$.

First let us assume that $\ann_R M \cap (m-m^2) \ne \phi$. Let
$x \ \in \ \ann_R M \cap (m-m^2)$ and let $\overline R = R/xR$.
Since $\pd_R M < \propto$ and $x  \in  m-m^2$, $\pd_{\overline R}M < \propto$.
By hypothesis M satisfies the order ideal conjecture as an $\overline R$-module.
Let $(F_\bullet, d_\bullet)$, $(P_\bullet, S_\bullet)$ denote minimal free
resolutions of M over R and $\overline R$ respectively; let $S_i = Syz^i_R (M)$
and $T_i = Syz^i_{\overline R} (M)$. Arguing as in the proof of part b) of the
theorem ((1), (2), (3) etc.) we construct a minimal set of generators
$\{xe_1, \dots, xe_i, \alpha_js\}$ such that $\{e_1, \dots, e_j\}$ form a part of a
basis of $F_0$, $\overline \alpha_j (= im \alpha_j  \in  T_1)$ form a
minimal set of generators of $T_1$. In order to show that for any minimal
generators $\beta$ of $S_1$, $\grade \ocal_{S_1} (\beta) \ge 1$, due to
inductive hypothesis, it is enough to take $\beta = xe - \Sigma \lambda_i \alpha_i,
\ e  \in \{e_1, \dots, e_i\}, \lambda_i \in  m$. Then $\beta_e$ = the e-th
co-ordinate of $\beta = x - \Sigma \lambda_i \alpha_{ie}$. Since
$\lambda_is  \in  m$, $\alpha_{ie} s  \in \ m$ and $x  \in  m-m^2$,
we have $\beta_e \ne 0$, $\ \beta_e  \in  m-m^2$ and hence by assumption
$\beta_e$ is a non-zero-divisor in R. Thus the conclusion follows for i = 1.

Now consider $i > 1$. Arguing as in part c) above we see from
$(1^\prime), (2^\prime), (3^\prime)$ etc. that it is enough to
consider a minimal generator $\beta$ of $S_i$ of the form
$\beta = \gamma - \Sigma \lambda_i \alpha_i$ where $\gamma = xe^\prime -
h_{i-2} h _{i-1} (e^\prime) = d_i (e), e = h_{i-1} (e^\prime)$, $e^\prime$ a minimal
generator of $F^\prime_{i-1}$ and ${\overline \alpha}_j \( =  im \alpha_j \ \in T_i \)$  form a minimal set of
generators of $T_i$. Similar arguments as in part c) show that grade of
$\ocal_{S_i} (\gamma) \ge i$ and by hypothesis, $\grade \ocal_{T_i} ({\overline \alpha}_j) \ge i$.
Let J denote the ideal generated by entries
of $\beta$. Recall that $F_{i-1} = F_{i-1}^\prime \oplus F_{i-1}^{\prime\prime}$,
$e^\prime$ is a minimal generator of $F_{i-1}^\prime, h_{i-2} d_{i-1}
(e^\prime) \ \in \ F_{i-1}^{\prime\prime}$. Let $y$ = the $e^\prime$-th
co-ordinate of $\beta = x - \Sigma \lambda_i \alpha_{ie^\prime}$; since
$\lambda_is \ \in \ m$ and $\alpha_{ie^\prime}s \ \in \ m$,
$y \ \in \ m-m^2$ and hence, by assumption, $y$ is a non-zero-divisor in R. Let
$\overline R = R/yR$, $\overline S_i = S_i/yS_i$ for $i \ge 1$; then
$\overline S_i = Syz_{\overline R}^{i-1} (\overline S_1)$. By hypothesis $J/yR$
has $\grade \ge (i-1) $ in $\overline R$. This implies that $\grade_R J \ge i$
and our proof is complete.

Now assume $\depth_R M > 0$. we can find $x  \in  m-m^2$
such that $x$ is a non-zero-divisor on M. Let $\overline R = R/xR$,
$\overline M = M/xM$. Since $\pd_{\overline M} = \pd_R M < \infty$, by hypothesis,
$\overline M$ satisfies the order ideal conjecture over $\overline R$.
And hence M satisfies the order ideal conjecture over R.

\enddemo

\proclaim {Corollary\enspace 1} Let $(R, m)$ be a local ring of mixed characteristic $p > 0$
such that $p$ is a non-zero-divisor in $R$. Let $I$ be an ideal of finite projective dimension over $R$.
If $p \in I$ or $p$ is a non-zero-divisor on $R/I$ then every minimal generator of $I$ is a
non-zero-divisor in $R$. In particular if $P$ is a prime ideal of finite projective dimension over $R$,
then every minimal generator of $P$ is a non-zero-divisor.
\endproclaim
\demo{Proof} If $p\in I$, the proof follows from part b) of the above theorem; if $p$ is a non-zero-divisor
on $R/I$, the result follows from the validity of the order ideal conjecture on $R/pR$ ([E-G2]).
\enddemo

\proclaim {Corollary\enspace 2}\enspace Let (R, m) be a regular local ring of dimension $n$ and assume that
the order ideal conjecture is valid for regular local rings of dimension $(n-1)$. If $M$ is a finitely
generated $R$-module such that either $M$ is annihilated by a regular parameter or $\depth_R M > 0$,
then $M$ satisfies the order ideal conjecture.
\endproclaim
Proof follows from part d) of the above theorem.

\enddocument